\documentclass[11pt,a4paper]{amsart}

\usepackage{latexsym}
\usepackage{fancyhdr}
\usepackage{amsmath,amsthm,amsfonts,amsbsy,amsgen,amscd,mathrsfs}
\usepackage{bbm}
\usepackage{euscript}
\usepackage{graphicx}
\usepackage[colorlinks=true]{hyperref}

\edef\qedrestoreat{\noexpand\catcode\lq\noexpand\@=\the\catcode\lq\@}
\catcode\lq\@=11

\ifx\protect\undefined\let\protect\relax\fi

\def\qed{\protect\@qed{$\qedsymbol$}}
\def\pushright{\protect\@pushright}

\def\QED{\protect\@qed{{\rm Q.E.D.}}}

\def\QEI{\protect\@qed{{\rm Q.E.I.}}}

\def\Proof{\protect\@Proof}\def\endProof{\protect\@endProof}%

\def\Proofof#1{\protect\@Proofof{#1}}\def\endProofof{\protect\@endProofof}%

\def\qedsymbol{\raisebox{-.2ex}{$\Box$}}

\def\TheWordProof{\sc Proof.}
\def\TheWordProofof#1{\sc Proof of #1.}

\def\ProofFont{}

\newif\ifAutoQED\AutoQEDfalse

\newif\ifNumberResults
\ifx\theorem@style\undefined
   \NumberResultstrue
   
\fi

\def\parag@pushright#1{{
    \parfillskip=0pt            
    \widowpenalty=10000         
    \displaywidowpenalty=10000  
    \finalhyphendemerits=0      
    \hbox@pushright             
    #1
    \par}}

\def\hbox@pushright{
    \unskip                     
    \nobreak                    
    \hfil                       
    \penalty50                  
    \hskip.2em                  
    \null                       
    \hfill                      
}%

\newif\if@qed\@qedfalse
\def\save@set@qed{\let\saved@ifqed\if@qed\global\@qedtrue}%
\def\restore@qed{\global\let\if@qed\saved@ifqed}

\def\@Proof{%
   \par\removelastskip\bigskip\penalty100
   \save@set@qed
   \noindent\ProofFont{\TheWordProof\enskip}%
}%
\def\@Proofof#1{%
   \par\removelastskip\bigskip\penalty100
   \save@set@qed
   \noindent\ProofFont{\TheWordProofof{#1}\enskip}%
}%
\def\@endProof{%
   \qed\restore@qed
   \penalty-100 \medskip
}
\def\@endProofof{%
   \qed\restore@qed
   \penalty-100 \medskip
}
\def\@qed#1{%
\if@qed                                 
     \global\@qedfalse
        \ifmmode\ifinner\pushright{#1}
        \else\eqno{\qedsymbol}\fi
        \else\pushright{#1}\fi%
\else\ifhmode\ifinner\else\par\fi\fi
\fi}
\def\@pushright#1{%
  {\ifvmode                       
       \null\hfill{#1}\par        
  \else\ifmmode\maths@pushright{\hbox{#1}}
       \else\ifinner\hbox@pushright{#1}
            \else\parag@pushright{#1}
  \fi  \fi  \fi
}}%
\def\maths@pushright#1{{%
  \ifinner
     \hbox@pushright{#1}%
  \else
     \eqno#1
     \def\]{$$\ignorespaces}
  \fi
}}%
\theoremstyle{plain}
\newtheorem{lemma}{Lemma}[section]	
\newtheorem{theorem}{Theorem}[section]

\newtheorem{corollary}{Corollary}[section]

\theoremstyle{definition}
\newtheorem{definition}{Definition}[section]

\theoremstyle{remark}

\newtheorem{example}{Example}[section]

\def\R{\mathbb{R}}

\def\Gr{\mathbb{G}}

\def\Prob{\mathop{\mathbf{P}}}

\DeclareMathOperator{\dist}{d}

\newcommand{\e}{\varepsilon}

\renewcommand{\>}{\rangle}

\def\Oh{{\mathcal O}}

\newcommand{\normal}{N}

\newcommand{\mdeg}{\operatorname{mdeg}}
\newcommand{\binomgr}[2]{\left[\begin{array}{c} #1\\#2\end{array}\right]}

\def\algorithm{\begin{center}
               \begin{minipage}{6in}
               \begin{tabbing}
               \marks}
\def\falgorithm{\end{tabbing}
                \end{minipage}
                \end{center}}
\def\marks{nn\= nn\= nn\= nn\= nn\= nn\= nn\= \kill}

\def\II{{\mathrm{I\kern-0.5pt I}}}

\def\vol{{\mathsf{vol}}}

\def\Prob{\mathop{\mathbf P}}

\newcommand{\binomial}[2]{\ensuremath{\left(
\begin{array}{c} #1 \\ #2 \end{array} \right)}}

\def\dist{{\sf dist}}

\def\vol{{\mathsf{vol}}}

\begin{document}

\title[On the volume of tubes]{On the Volume of Tubular Neighborhoods of Real Algebraic Varieties}
\author{Martin Lotz}

\address{School of Mathematics, The University of Manchester, Alan Turing Building,
Oxford Road, Manchester, M139PL, United Kingdom.}
\thanks{Research supported by Leverhulme Trust
    grant R41617 and a Seggie Brown Fellowship of the University of
    Edinburgh}


\begin{abstract}
The problem of determining the volume of a tubular neighborhood has a long and rich history. Bounds on the volume of neighborhoods of algebraic sets
have turned out to play an important role in the probabilistic analysis of condition numbers in numerical analysis.
We present a self-contained derivation of bounds on the probability that a random point, chosen uniformly from a ball, lies within a given distance of a real algebraic variety of any codimension.
The bounds are given in terms of the degrees of the defining polynomials, and contain as special case an unpublished result by Ocneanu.
\end{abstract}

\maketitle

\section{Introduction}
The purpose of these notes is to derive a bound on the volume of a tubular neighborhood of a real algebraic variety in terms of the degrees of the defining polynomials.
The problem is stated in probabilistic terms, namely, as the probability that a random point, uniformly distributed in a ball, falls within a certain neighborhood of the variety.

\begin{theorem}\label{th:main}
Let $V$ be the zero-set of multivariate polynomials $f_1,\dots,f_s$ in $\R^n$ of degree at most $D$.
Assume $V$ is a complete intersection of dimension $m=n-s$. Let $x$ be uniformly distributed in a ball $B^n(p,\sigma)$ of radius $\sigma$ around $p\in \R^n$. Then
\begin{equation*}
  \Prob\{\dist(x,V)\leq \e\}\leq
  4 \ \sum_{i=0}^{m}\binomial{n}{s+i} \ 
\left(\frac{2D\e}{\sigma}\right)^{s+i}\left(1+\frac{\e}{\sigma}\right)^{m-i}.
\end{equation*}
If the polynomials $f_1,\dots,f_s$ are homogeneous and $p=0$, then
\begin{equation*}
\Prob\{\dist(x,V)\leq \e\}\leq
 2 \ \sum_{i=0}^{m}\binomial{n}{s+i} \ \left(\frac{2D\e}{\sigma}\right)^{s+i}.
\end{equation*}
\end{theorem}

The second of the stated equations is commonly attributed to A.~Ocneanu~\cite[Theorem 4.3]{demm:88}, though a proof has not been published so far 
and does not seem available. From the proof of Theorem~\ref{th:main} we also get the following corollary, conjectured by J.~Demmel~\cite[(4.15)]{demm:88}.

\begin{corollary}\label{co:asympt}
For compact $V$ and small enough $\e$ we have
\begin{equation*}
\Prob\{\dist(x,V)\leq \e\} = \vol_{n-s}(V) \cdot \e^s\cdot \frac{n\Gamma(n/2)}{\pi^{(n-s)/2}s\Gamma(s/2)}+o(\e^{s}).
\end{equation*}
\end{corollary}

Theorem~\ref{th:main} can be adapted to a spherical setting without too much difficulty, thus generalizing the results of~\cite{bucl:08} to higher codimension, but for
the sake of brevity such a generalization is omitted in these notes.

\subsection{History and applications}
In 1840, J. Steiner~\cite{stei:40} showed that volume of an $\e$-neighborhood of a convex body in $\R^3$ could be written
as a quadratic polynomial in $\e$. This result has become a staple of integral geometry and was the starting point of a myriad of generalizations in multiple directions.
One such generalization is a celebrated result by H. Weyl~\cite{weyl:39}, who showed that for $\e$ small enough, the volume of an $\e$-neighborhood around a compact Riemannian submanifold of $\R^n$ is given by 
a polynomial whose degree is the dimension of the manifold.
Weyl's tube formula became an important ingredient in Allendoerfer and Weil's proof of the Gauss-Bonnet Theorem for hypersurfaces. For more on Weyl's tube formula and its ramification, see~\cite{gray:04}.
Bounds on the volume of tubes around real varieties in terms of degrees have previously been given by R.~Wongkew~\cite{wong:93}, although without explicit constants. 
Tube formulae came into the radar of numerical analysis through the work of S. Smale~\cite{smal:81}, E. Kostlan~\cite{kost:85}, J. Renegar~\cite{rene:87}, and J. Demmel~\cite{demm:88}, among others, who were interested in the
probabilistic analysis of condition numbers. It has been observed (see, e.g.,~\cite{kaha:72,demm:87} and the references there) that the condition number of many numerical computation problems
can by bounded by the inverse distance to a set of ill-posed inputs.
In particular, if one can describe the set of ill-posed inputs as a subset of an algebraic variety, then a bound on the relative volume of its neighborhood in terms of the degree of the variety directly translates
into a result on the probability distribution of condition numbers. The results of Demmel~\cite{demm:88} have been partially extended to the setting of {\em smoothed analysis} on the sphere in~\cite{bucl:08},
by studying tubular neighborhoods of hypersurfaces intersected with spherical caps. 
For a comprehensive survey of these ideas we refer to~\cite{buer:10}.
Recently, a consequence of the degree bound derived in this article
has been used in the study of embeddings of simplicial complexes into Euclidean space~\cite[Prop 3.10]{grgu:12}. 
Other notable fields in which tube formulae have been used extensively include statistics~\cite{adta:07}
and the probabilistic analysis of convex optimization~\cite{ambu:11,almt:13,mctr:13}.
The main purpose of the current article is to fill a gap in the literature by making available a complete and rigorous derivation of the real degree bounds used in~\cite{demm:88}.

\subsection{Main ideas}
The proof of Theorem~\ref{th:main} is based on three main ingredients: Weyl's tube formula, an integral-geometric kinematic formula, and B\'ezout-type 
bounds on the degree of Gauss maps. In what follows, let $V$ be a complete intersection of dimension $m=n-s$.
First, based on Weyl's tube formula, a bound is derived in terms of {\em integrals of absolute curvature}:
\begin{equation*}
\vol_n \ T(V,\e) \leq \e^s \sum_{i=0}^{m} \frac{1}{s+i} \ |K_{i}|(V) \ \e^{i}.
\end{equation*}
The highest order term $|K_m|(V)$ is intimately related to the {\em generalized Gauss map} of $V$, and can in fact be expressed in terms of the {\em degree} of this map. 
Using standard B\'ezout-type arguments it is possible to bound the degree of the Gauss map in terms of the degrees of the defining polynomials.
The lower-order invariants $|K_i|(V)$ can then be related to the highest order invariants $|K_i|(V\cap L)$ of an intersection with a random linear subspace by means of {\em Crofton's Formula} from integral geometry:
\begin{equation*}
  |K_i|(V)\leq 2\binomgr{n}{s+i}\int_{L\in \mathcal{E}_{s+i}^n}|K_i|(V\cap L)
  \ d\lambda_{s+i}^n.
\end{equation*}
where $\mathcal{E}_{s+i}$ denotes the space of $(s+i)$-dimensional subspaces with suitable measure.
One can the apply the degree bounds in lower dimension. Obviously, some care has to be taken
when implementing these ideas in detail.

\subsection{Outline}
Section~\ref{se:prel} gives a review of the necessary concepts of Riemannian geometry in Euclidean space.
In Section~\ref{se:tubes}, Weyl's tube formula and results from integral geometry are presented in a slightly generalized form to suit our purposes.
At the beginning of Section~\ref{se:degree}, the tube formula is reformulated in terms of the degrees of a generalized Gauss map.
Up to this point, everything is based on compact Riemannian manifolds.
Systems of polynomial equations enter when bounding the degrees of the generalized Gauss map,
leading to the proof of Theorem~\ref{th:main}.
The appendix is devoted to a complete proof of Weyl's tube formula in Euclidean space.

\subsection{Notation and terminology}\label{se:notation}
We write $B^n(p,\sigma)$ for the solid closed ball in $\R^n$ with
center $p$ and radius $\sigma>0$,
and $S^{n-1}(p,\sigma)$ for its boundary, and set $S^{n-1}:=S^{n-1}(0,1)$ and $B^n:=B^n(0,1)$.
We write $\vol_n\ M$ for the $n$-dimensional
Lebesgue-measure of a measurable set $M\subseteq \R^n$, and often drop
the subscript an simply write $\vol \ M$. For an $m$-dimensional Riemannian manifold
$M$, when we write $\vol_m \ M = \vol \ M$ we mean $\int_M\omega_M$,
with $\omega_M$ the volume form associated to the Riemannian
structure (see Section~\ref{sec:integration}). Whenever we say manifold, we mean smooth manifold.

Throughout this paper we denote by
$\Oh_{n-1}:=2\pi^{n/2}/\Gamma\left(\frac{n}{2}\right)$ the
$(n-1)$-dimensional volume of the unit sphere $S^{n-1}$ in $\R^{n}$, and
$\omega_n:=\Oh_{n-1}/n$ the $n$-dimensional volume of the solid unit
ball in $\R^n$. The {\em flag coefficients} are defined as
\begin{equation}\label{eq:bingr}
  \binomgr{n}{k}:=\binomial{n}{k}\frac{\omega_n}{\omega_k \ \omega_{n-k}}
\end{equation}
for $n\geq 0$ and $k\geq 0$. They appear naturally in the study of invariant
measures on Grassmannians~\cite{klro:97}.

\section{Preliminaries}\label{se:prel}
We assume familiarity with the basic notions of Riemannian geometry, as described for example in~\cite{doca:92}.
The purpose of most of this section is to introduce notation and terminology.

\subsection{Riemannian manifolds in $\R^n$}\label{sse:riemannian}
Given a Riemannian manifold $M$ of dimension $m$, we denote by $TM$ its tangent bundle, by $C(M)$ the ring of smooth functions on $M$,
and by $\mathcal{X}(M)$ the $C(M)$-module of tangent vector fields on
$M$. For $p\in M$ we write $T_pM$ for the tangent space at $p$. 
If $v\in T_p\R^n$ and $f\in C(\R^n)$, then $v(f)$ denotes the directional derivative of $f$ in direction $v$ at $p$.

In this article we are only concerned with submanifolds $M$ of
Euclidean space $\R^n$. As such, each $T_pM$ can be identified with a
subspace of $T_p\R^n\cong \R^n$ in the obvious manner. 
Let $\normal M := \{ (p,v)\in T\R^n \mid p\in M, v \perp T_pM\}$
be the normal bundle to $M$ in $\R^n$ and denote by $\normal_pM$ the
fiber of $\normal M$ over $p\in M$, i.e., the normal space to $M$ at $p$ in $\R^n$.

Let $Y\in \mathcal{X}(\R^n)$ be a smooth vector field. For $v\in T_p\R^n$ denote by $\nabla_vY:=v(Y)$ the covariant derivative of
$Y$ along $v$ at $p$. The covariant derivative satisfies $v(\langle Y,Z\rangle)
=\langle\nabla_vY,Z\rangle+\langle Y,\nabla_vZ\rangle$. In particular,
for orthogonal fields $Y$ and $Z$ we have $\langle
\nabla_vY,Z\rangle=-\langle Y,\nabla_vZ\rangle$. 
For $v\in \normal_p M$ and $X,Y\in \mathcal{X}(M)$, the second
fundamental form $S_v(X,Y)$ of $X$ and $Y$ along $v$ is the symmetric,
bilinear map $T_pM\times T_pM\rightarrow \R$ defined by
\begin{equation*}
  S_v(X,Y):=\langle \nabla_{X(p)}Y,v\rangle,
\end{equation*} 
where we assume the vector fields $X,Y$ to be extended to a
neighborhood of $M$ in $\R^n$ for this definition to make sense.
Given a normal vector field $Z$ on $M$ we have $S_{Z(p)}(X,Y)=-\langle
Y,\nabla_{X(p)}Z\rangle$ (since $X,Y$ are orthogonal to $Z$).
Given an orthonormal frame field $(E_1,\dots,E_m)$ on $U\subset M$, we will on
occasion use the matrix $S(Z)$ with entries in $C(M)$ that represents this bilinear form with respect to that frame field. Its values at
$p\in U$ are given by the entries
\begin{equation}\label{eq:second}
  S_{ij}(Z)(p)=S_{Z(p)}(E_i,E_j)=\langle \nabla_{E_i(p)}E_j,Z\rangle=-\langle E_j,\nabla_{E_i(p)}Z\rangle.
\end{equation}
Note that we can also talk about $S(v)$ for fixed $v\in \normal_pM$. Then we have
$S(v)=S(Z)(p)$ for any normal vector field such that $Z(p)=v$.

\subsubsection{A note on integration and orientation}\label{sec:integration}
Given a Riemannian manifold $M\subseteq \R^n$, we denote by $\omega_M$ the natural volume form on $M$ associated to the Riemannian metric. 
Thus if $U\subseteq M$ is an oriented coordinate neighborhood and $x^1,\dots,x^m\colon U\rightarrow \R^m$
are orthonormal coordinates (so that the tangent vectors $\partial/\partial x^1,\dots,\partial/\partial x^m$ form a positively oriented,
orthonormal basis at each $p\in U$), then $\omega_M=dx^1\wedge \cdots \wedge dx^m$ on $U$. All volume forms are densities (unsigned forms), though
we will occasionally locally represent them as differential forms in an oriented coordinate neighborhood $U\subseteq M$ without always stating this explicitly.
Given a map $f$ from a manifold of the same dimension to $M$, $f^*\omega_M$ denotes the pull-back volume form.

\subsubsection{Curvature Invariants}\label{sec:curvinv}
In this section we introduce the curvature invariants $K_0(M),\dots,K_m(M)$ associated to a compact Riemannian manifold $M$ in terms
of the second fundamental form. These invariants are key components in Weyl's formula (Section~\ref{sec:weyl}) for the volume of tubes around $M$.

Let $M$ be an $m$-dimensional compact Riemannian manifold, $U'\subseteq \R^n$ open and $U=U'\cap M$.
Let $(E_1,\dots,E_n)$ be an orthonormal frame field on $U'$, such that $E_1(p),\dots,E_m(p)$ form an oriented basis of
$T_pM$ for all $p\in U$. For $v\in \normal_pM$ let $S(v)$ denote the $m\times m$ matrix of the
second fundamental form at $p$ along $v$ with respect to the frame field, as defined in~(\ref{eq:second}).

For $0\leq i\leq m$ let $\psi_i\colon \normal_p M\rightarrow \R$ be the homogeneous polynomial of
degree $i$ defined by
\begin{equation*}
  \det(\mathrm{Id}-tS(v))=\sum_{i=0}^m t^i \psi_i(v).
\end{equation*}
Note that the $\psi_i(v)$ are, up to sign, the
coefficients of the characteristic polynomial of $S(v)$. More
precisely, we have
\begin{equation*}
  \psi_i(v)=(-1)^i\sigma_i(\kappa_1(v),\dots,\kappa_m(v)),
\end{equation*}
where the $\kappa_i(v)$ are the eigenvalues of $S(v)$, i.e., the
principal curvatures, and $\sigma_i$ denotes the $i$-th elementary
symmetric function.
In particular, $\psi_m(v)=(-1)^m\det S(v)$. These quantities are, up to orientation, independent of the particular orthonormal
frame used to define the matrix $S(v)$.

For $p\in M$, set $S(\normal_pM):=\{v\in \normal_pM \mid \|v\|=1\}$ and denote by $S(NM)$ the corresponding normal sphere bundle. At a point $p\in M$ define
\begin{equation*}
  I_i(p)=\int_{v\in S(\normal_pM)} \psi_i(v)\; \omega_{S(\normal_pM)}.
\end{equation*}
The $I_i(p)$ are polynomial invariants of the second fundamental form in the sense of~\cite{howa:93}.
If $v=\sum_{j=1}^su^jE_{m+j}(p)$, $s=n-m$, then the $I_i(p)$ are integrals over
all $u\in S^{s-1}$ of homogeneous polynomials of degree
$i$ in $u^1,\dots,u^s$. From this it follows that $I_i(p)=0$ for $i$ odd. 

The {\em integrals of curvature} are defined as
\begin{equation}\label{eq:defki}
  K_i(M):=\int_M I_i(p)\; \omega_M=\int_{S(\normal M)}\psi_i(v)\;
  \omega_{S(\normal M)}.
\end{equation}
It is easy to see that $K_0(M)=\Oh_{s-1}\vol_m \ M$. Less trivial is
the fact that $K_m(M)=\Oh_{n-1}\ \chi(M)$, where $\chi(M)$ is the Euler
characteristic of $M$. This is a consequence of the generalized
Gauss-Bonnet Theorem (see~\cite{gray:04} for a discussion of this
result and its relation to Weyl's tube formula).

The {\em integrals of absolute curvature}, suggested by Peter B\"urgisser~\cite{bucl:08}, are defined as
\begin{equation}\label{eq:defabski}
  |K_i|(M):=\int_{S(\normal M)}|\psi_i(v)|\; \omega_{S(\normal M)}.
\end{equation}
These are important for extending Weyl's tube formula to and inequality for the volume of $\e$-tubes 
for arbitrary $\e$. Clearly, the definition is also valid for an open subset $U\subset M$, or an open subset of $M\backslash \partial M$ if
$M$ is a compact Riemannian manifold with boundary.

%
%
%
\subsubsection{The degree}
Let $f\colon M\rightarrow P$ by a smooth map of compact Riemannian manifolds.
By Sard's Theorem~\cite[\S 2]{miln:97} almost all $q\in P$ are regular values.
The preimage $f^{-1}(q)$ is either empty or a
finite set with locally constant cardinality~\cite{miln:97} as $q$ varies among regular values. 

For measurable $h\colon P\rightarrow \R$ we have
\begin{equation}\label{eq:coarea}
  \int_{p\in M}h\circ f(p) \ f^*\omega_P = \int_{q\in P} h(q)\ \#f^{-1}(q) \ \omega_P,
\end{equation}
where $\#f^{-1}(q)$ denotes the cardinality of the preimage of $q$.
Recall (Section~\ref{sec:integration}) that we are dealing with unsigned forms, i.e., $f^*\omega_P = |\det(D\varphi)|\omega_m$, otherwise we would have to count the
points in the fiber with signs.

We define the maximum degree of $f$ to be the maximum cardinality of the preimage of a regular value under $f$:
\begin{equation*}
  \mdeg f := \max_{q\in \mathrm{reg }P} \#f^{-1}(q).
\end{equation*}
With this definition we have
\begin{equation}\label{eq:degree}
  \int_{M}\ f^*\omega_P\leq \mdeg f \ \int_{P}\ \omega_P.
\end{equation}

This notion of degree differs from the usual one from differential topology (see~\cite[\S 5]{miln:97}), 
which takes into account orientation. 


\subsubsection{Transversality}
The intersection of two manifolds $M$ and $P$  in $\R^n$ of dimension $m,\ell$ with
$m+\ell\geq n$ is called {\em transversal} at $p\in M\cap P$, if $\dim
T_pM\cap T_pP=m+\ell-n$. The intersection is called transversal if it
is transversal at every $p\in M\cap P$. In that case, $M\cap P$ is a
smooth $(m+\ell-n)$-dimensional manifold. 

Recall that $B^n(p,\sigma)$ denotes the closed ball of radius $\sigma$
around $p$ in $\R^n$, and $S^{n-1}(p,\sigma)=\partial B^n(p,\sigma)$
is its boundary. The following lemma is a standard application of Sard's Lemma,
the proof is omitted. By ``almost all'' we mean ``up to a set of measure zero''.

\begin{lemma}\label{le:transverse}
Let $M$ be a Riemannian manifold. For almost all $\sigma>0$ the
intersection $B^n(p,\sigma)\cap M$ is a Riemannian manifold with boundary.
In particular, $S^{n-1}(p,\sigma)\cap M$ is a smooth Riemannian manifold of
codimension one in $M$.
\end{lemma}


\section{Geometry of tubes and integral geometry}\label{se:tubes}

\subsection{Weyl's tube formula}\label{sec:weyl}
References for the content of this section are~\cite{weyl:39,gray:04}.
Let $M\subseteq \R^n$ be a Riemannian submanifold of dimension
$m< n$, possibly with boundary, and denote by $s:=n-m$ the codimension of $M$ in $\R^n$.
The (closed) {\em tube} of radius $\e$ around $M$ in $\R^n$ is defined to be the set
\begin{equation}\label{eq:tube}
  T(M,\e) := \left\{p\in \R^n \bigg| \begin{array}{ll} \exists \text{
    line of length }\leq \e \text{ from } p \\ \text{ meeting } M
  \text{ orthogonally}\end{array}\right\}.
\end{equation}
For compact manifolds this coincides with the $\e$-neighborhood of $M$
in $\R^n$, though in general this need not be the case.

\begin{figure}[h]
\begin{center}
\includegraphics[width=0.8\textwidth]{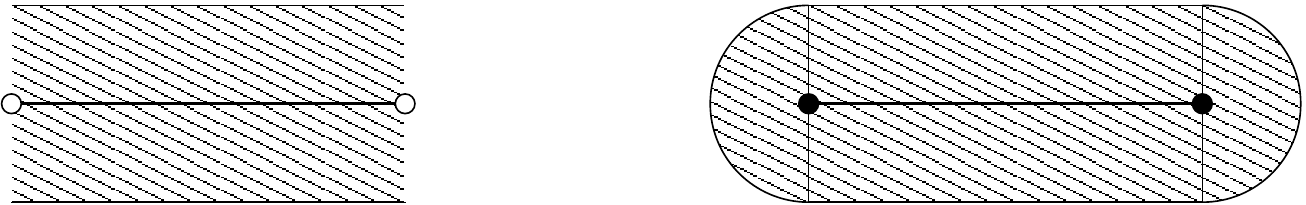}
\end{center}
\caption{Tube around and open [left] and closed [right] interval.}
\end{figure}

In his influential paper~\cite{weyl:39}, Weyl derived the expression
\begin{equation*}
\vol \ T(M,\e) = \Oh_{s-1} \ \e^s \sum_{\stackrel{i=0}{i \text{ even}}}^{m} 
\frac{(i-1)(i-3)\cdots 1}{(s+i)(s+i-2)\cdots s} \ \mu_{i}(M)\ \e^{i}
\end{equation*}
for the volume of a tube of radius $\e$ around $M$, provided $\e$ is
small enough. The $\mu_i(M)$ are the {\em curvature invariants}
of $M$. The deeper part of Weyl's work consists of showing that these
invariants are intrinsic, that is, they only depend on the curvature
tensor of $M$ and not on the particular embedding of $M$ in $\R^n$. We
will not need this feature here, however, and will be happy with
expressing these invariants
in terms of the second fundamental form.

The $\mu_i(M)$ are just a different normalization of the invariants $K_i(M)$ introduced in Section~\ref{sec:curvinv}, namely
\begin{equation}\label{eq:muk}
  K_i(M)=\Oh_{s-1}\frac{(i-1)(i-3)\cdots 1}{(s+i-2)(s+i-4)\cdots s} \mu_i(M).
\end{equation}
for $i$ even. Note that Corollary~\ref{co:asympt} follows immediately from the Weyl's tube formula, using that $K_0(M)=\Oh_{s-1}\vol_m\ M$.

Note that the $K_i(M)$ are no longer independent of the embedding,
since the codimension enters the formula. 

We will need a slight variation of Weyl's tube formula that works for 
arbitrary $\e$. 

\begin{theorem}\label{thm:tubeineq}
Let $M\subseteq \R^n$ be an oriented, compact, $m$-dimensional Riemannian
manifold, possibly with boundary, and assume $s:=n-m>0$. Let $U\subseteq M\backslash \partial M$ be an open subset of $M$. Then for all $\e>0$ we have
\begin{equation}\label{eq:weyl}
\vol \ T(U,\e) \leq \e^s \sum_{i=0}^{m} \frac{1}{s+i} \ |K_{i}|(U) \ \e^{i}.
\end{equation}
\end{theorem}

The proof, given in the appendix, is along the lines of ~\cite{weyl:39,gray:04}.

\begin{example}
Let $M=S^m$ be the $m$-dimensional unit sphere in $\R^n$. Along the lines of the
proof of the tube formula~\ref{thm:tubeineq} we can derive
\begin{equation*}
\vol \ T(S^m,\e) = 2 \ \omega_{m+1}\ \e^s \sum_{i=0}^{m} \frac{\omega_{s+i}}{\omega_{i+1}}\binomial{m+1}{i+1}\ \e^{i}.
\end{equation*}
for small $\e$ (recall from Section~\ref{se:notation} the definition of $\omega_n$ and $\Oh_n$). From this we get 
\begin{equation*}
K_i(S^m)=\frac{2\Oh_m\Oh_{s+i-1}}{\Oh_i}\binomial{m}{i}
\end{equation*}
for $i$ even.
Note that $K_0(S^m)=\Oh_m\Oh_{n-m-1}$ and
that $K_m(S^m)=2\Oh_{n-1}$ for $m$ even and $K_m(S^m)=0$ for $m$ odd,
in accordance with the Euler characteristic for spheres. Some
special cases for the volume of tubes:
\begin{enumerate}
\item Setting $m=n-1$ we get $\vol \  T(S^{n-1},\e)=\omega_n \ [(1+\e)^n-(1-\e)^n]$,
as was to be expected.
\item For $m=0$ we have $\vol \ T(S^0,\e)=2\ \e^n\ \omega_n$.
\item For $m=1$, $n=2$ we get the volume of the torus $\vol
  \ T(S^1,\e)=2\pi^2\e^2$.
\end{enumerate}
\end{example}

\subsection{Integral geometry}
In order to obtain upper bounds for the integrals of absolute curvature $|K_i|(M)$, we will first derive bounds for
$|K_m|(M)$ using the generalized Gauss map, and the case where $0\leq i<m$ is then handled by relating the $i$-th curvature invariants
$K_i(M)$ to the curvature invariants $K_i(M\cap L)$ of the
intersection of $M$ with a random affine space of dimension
$s+i$. Formulae relating invariant measures of a set to its intersection with
random affine spaces are known by the name of {\em Crofton formulae} and play a central role in integral geometry.
For an introduction to integral geometry and geometric probability
we refer to~\cite{klro:97, shwe:08}. The version of
Crofton's formula involving Weyl's curvature invariants is due to
Chern~\cite{cher:66} and Federer~\cite{fede:59}, see
also~\cite[15.95b]{sant:76}.

Let $\mathcal{E}_k^{n}$ be the set of $k$-dimensional affine spaces
in $\R^n$ and $\Gr(n,k)$ the Grassmannian of $k$-dimensional linear
subspaces of $\R^n$. We can identify $\mathcal{E}_k^n$ with the
subset of those $(V,p)\in \Gr(n,k)\times \R^n$ such that $p\perp V$,
the one-to-one correspondence $s$ being given by $s(V,p)=p+V$~\cite[Chapter 6]{klro:97}. 

Let $\nu_k^n$ denote the $O(n)$-invariant measure on $\Gr(n,k)$
induced by the identification $\Gr(n,k)=O(n)/O(k)\times O(n-k)$,
normalized such that
\begin{equation*}
  \nu_k^n(\Gr(n,k))=\frac{\Oh_{n-1}\cdots \Oh_{n-k}}{\Oh_{k-1}\cdots \Oh_0},
\end{equation*}
see also ~\cite[3.2]{bucl:09} for a discussion.
The product measure $\nu_{k}^n\times \omega_{\R^n}$ gives rise to an invariant
measure $\overline{\lambda}_k^n$ on $\mathcal{E}_k^n$, defined by
\begin{equation*}
  \int_{V\in \Gr(n,k)}\left(\int_{p\in
    V^{\perp}}f\circ s(V,p)\ \omega_{V^{\perp}}\right)\ d\nu_k^n=\int_{L\in
    \mathcal{E}_k^n} f(L)\ d\overline{\lambda}_k^n.
\end{equation*}
for a measurable function $f$ on $\mathcal{E}_k^n$. In particular, setting
\begin{equation*}
f=\mathbf{1}_{B^n(p,\sigma)}=\begin{cases}
   1 & L\cap B^n(p,\sigma)\neq \emptyset\\
   0 & \text{ else }
\end{cases}
\end{equation*}
we get
\begin{align*}
\overline{\lambda}_{k}^n(\{L\in \mathcal{E}_{k}^n \mid L\cap
  B^n(p,\sigma)\neq\emptyset\})&=\int_{L\in \mathcal{E}_k^n}
  f(L)\ d\overline{\lambda}_k^n\\
&=\int_{V\in \Gr(n,k)}\left(\int_{p\in
    V^{\perp}}\mathbf{1}_{B^n(p,\sigma)}\ \omega_{V^{\perp}}\right)\ d\nu_k^n\\
  &=\omega_{n-k}\sigma^{n-k}\nu_k^n(\Gr(n,k)).
\end{align*}

In the following we use the renormalized measure $\lambda_k^n=\nu_k^n(\Gr(n,k))^{-1} \ \overline{\lambda}_k^n$, so that
\begin{equation*}
  \lambda_{k}^n(\{L\in \mathcal{E}_{k}^n \mid L\cap
  B^n\neq\emptyset\})=\omega_{n-k}.
\end{equation*}

The following theorem is merely a reformulation of~\cite[15.95b]{sant:76} with a different normalization of the measure,
and after simplifying the constants. Note also that with the
parameters chosen here, it makes no difference whether we formulate
this theorem with $\mu_i(M)$ or with $K_i(M)$, since $M\cap L$ has
generically the same codimension $s$ in $L$ as $M$ in $\R^n$. Recall
the definition~(\ref{eq:bingr}) of the flag coefficients in Section~\ref{se:notation}.

\begin{theorem}(Crofton's Theorem)\label{thm:crofton}
\begin{equation}\label{eq:crofton}
  K_i(M)=\binomgr{n}{s+i}\int_{L\in \mathcal{E}_{s+i}^n}K_i(M\cap L)
  \ d\lambda_{s+i}^n.
\end{equation}
\end{theorem}

Crofton's Theorem leads to a bound on integrals of absolute curvature.

\begin{theorem}\label{thm:croftonabs}
Let $M$ be a compact Riemannian submanifold of $\R^n$ of dimension $m<n$, and let $i\leq m$. Then
\begin{equation}\label{eq:croftonabs}
  |K_i|(M)\leq 2\binomgr{n}{s+i}\int_{L\in \mathcal{E}_{s+i}^n}|K_i|(M\cap L) \ d\lambda_{s+i}^n.
\end{equation}
\end{theorem}

\begin{Proof}
Let $M_+$ and $M_-$ denote the parts of $M$ on which $I_i(p)$ is positive and negative, respectively. Then $|K_i|(M)=|K_i(M_+)|+|K_i(M_-)|$.
\end{Proof}

\section{Degree bounds}\label{se:degree}

\subsection{Degree of the Gauss map}
In this section we interpret the expected value of the highest curvature invariant as the degree of a
generalized Gauss map. Let $S(\normal M)$ denote the normal sphere
bundle over $M$. Note that $S(\normal M)$ has codimension one in $\R^n$.

\begin{definition}
Let $M\subseteq \R^n$ be a compact $m$-dimensional Riemannian manifold.
The {\em generalized Gauss map} of $M$ is defined as
\begin{equation*}
  \gamma\colon S(\normal M)\rightarrow S^{n-1}, \quad (p,v)\mapsto v.
\end{equation*}
\end{definition}


%

The generalized Gauss map on a compact manifold can be shown to be surjective.
Note that for almost all $w\in S^{n-1}$, the map $h(p,v)=\<v,w\>$ is a Morse function with non-degenerate critical points those
$(p,v)$ such that $v=w$.

Recall now the definition~(\ref{eq:degree}) of the degree of a map. 

\begin{lemma}\label{le:curvedegree}
Let $M\subseteq \R^n$ be a compact Riemannian manifold of dimension $m$. Then
\begin{equation}\label{eq:gauss}
  |K_m|(M)=\int_{v\in S^{n-1}} \# \gamma^{-1}(v) \ \omega_{S^{n-1}} \leq \Oh_{n-1} \ \mdeg \gamma .
\end{equation}
\end{lemma}

\begin{Proof}
We need to show that 
\begin{equation}\label{eq:pullback}
  \gamma^*\omega_{S^{n-1}}=|\det S(v)| \ \omega_{S(\normal M)}
\end{equation}
on $M$. Once this is shown, the claim of the lemma follows from the definition
of the $|K_i|$ (\ref{eq:defki}), namely,
\begin{equation*}
  |K_m|(M)=\int_{(p,u)\in S(\normal M)} |\det S(u)|\ \omega_{S(\normal M)}.
\end{equation*}

Let $x^1,\dots,x^m\colon U\rightarrow \R^m$ be
orthonormal coordinates on an open set $U\subset M$. Let $(E_1,\dots,E_n)$ be an
orthonormal frame field defined in a neighborhood of $U$ in $\R^n$,
such that on $U$ we have $E_i=\partial/\partial x^i$ for $1\leq i\leq m$. 
The frame field $E_1,\dots,E_n$ gives a local trivialization of the sphere bundle
\begin{align*}
  U\times S^{s-1}&\rightarrow S(\normal M)\\
  (p,u)&\mapsto \left(p,\sum_{i=1}^{s}u^i E_{m+i}\right).
\end{align*}
An orthonormal coordinate system $y^1,\dots,y^{s-1}$ for $S^{s-1}$ thus gives rise to orthonormal coordinates
$x^1,\dots,x^m,y^1,\dots,y^{s-1}$ on $S(\normal M)$. 
With $\omega_M=dx^1\wedge \cdots \wedge dx^{m}$ and $dy=dy^1\wedge \cdots \wedge dy^{s-1}$ we have
\begin{equation}
  \omega_{S(\normal M)}=\omega_M\wedge dy.
\end{equation}
Similarly we have $\omega_{S^{n-1}}=E_1^*\wedge \cdots \wedge E_m^*\wedge dy^1\wedge \cdots \wedge dy^{s-1}$. 
Let $\phi$ be such that $\gamma^*\omega_{S^{n-1}}=\phi(p,v) \ \omega_{S(\normal M)}$ as differential form.
Then
\begin{align*}
  \phi(p,v)&=\gamma^*\omega_{S^{n-1}}\left(
\frac{\partial}{\partial x^1},\dots,\frac{\partial}{\partial
  x^m},\frac{\partial}{\partial y^1},\dots,\frac{\partial}{\partial
  y^{s-1}}\right)\\
&=\omega_{S^{n-1}}\left(\gamma_*\frac{\partial}{\partial x^1},\dots,\gamma_*\frac{\partial}{\partial
    x^m},\gamma_*\frac{\partial}{\partial
    y^1},\dots,\gamma_*\frac{\partial}{\partial y^{s-1}}\right).
\end{align*}
Note that
\begin{equation*}
\gamma_*\frac{\partial}{\partial x^i}=\sum_{\ell=1}^s u^{\ell}\frac{\partial}{\partial x^i}E_{m+\ell}(p), \quad
\gamma_*\frac{\partial}{\partial y^j}=\sum_{\ell=1}^s \frac{\partial}{\partial y^j}u^{\ell} \ E_{m+\ell}(p),
\end{equation*}
from which we obtain
\begin{equation*}
 \phi(p,v) = \omega_{M}\left(\frac{\partial}{\partial x^1} \gamma,\dots,\frac{\partial}{\partial x^m} \gamma\right) \cdot
  dy\left(\frac{\partial}{\partial y^1} \gamma,\dots,\frac{\partial}{\partial y^{s-1}}\gamma\right).
\end{equation*}

A direct calculation shows that
\begin{equation*}
  \langle\frac{\partial}{\partial x^i}\gamma,E_j\rangle =-S_{ij}(v).
\end{equation*}
From this it follows that 
\begin{equation*}
\omega_{M}\left(\frac{\partial}{\partial x^1}
\gamma,\dots,\frac{\partial}{\partial x^m} \gamma \right)=(-1)^m\det S(v).
\end{equation*}
Clearly
\begin{equation*}
dy\left(\frac{\partial}{\partial
    y^1} \gamma,\dots,\frac{\partial}{\partial y^{s-1}}
  \gamma \right)=1
\end{equation*}
from which the claim follows for $M$ without boundary. 
\end{Proof}

The statement of Lemma~(\ref{le:curvedegree}) also holds if $M$ is replaced by an open subset $U\subset M\backslash \partial M$, for a compact
manifold with boundary $M$. We omit the details.

For an affine subspace $L\in \mathcal{E}_{s+i}^n$ in general position,
the intersection $M\cap L$ is either empty or an $i$-dimensional
submanifold of $L\cong \R^{s+i}$. In the latter case we can define the
degree of $M$ with respect to $L$ as the degree of the Gauss
map of $M\cap L$ in $L$, that is,
\begin{equation*}
  \mdeg(M;L):= \mdeg \gamma|_{M\cap L}\leq \max_{v\in
    S^{s+i-1}}\#\gamma|_{M\cap L}^{-1}(v).
\end{equation*}  
Define the $i$-th degree $\mdeg_i(M)$ of $M$ to be the maximum of
$\mdeg(M;L)$ over all $L\in \mathcal{E}_{s+i}^n$ that intersect $M$:
\begin{equation*}
  \mdeg_i(M):=\sup_{L\in\mathcal{E}_{s+i}^{n}}\mdeg(M;L).
\end{equation*}

Before dealing with polynomial equations, we give a bound of the volume of an $\e$-tube around $M$ by a function
of the $i$-th degrees of $M$ and of $\e$. 

\begin{theorem}\label{thm:degbound}
Let $M$ be a Riemannian manifold of dimension $m$ in $\R^n$
and set $s=n-m$. Assume $M$ is contained in a ball $B^n(p,\sigma)$ of radius $\sigma$. Then for $\e>0$ we have
\begin{equation*}
\vol \ T(M,\e)\leq 2 \omega_n\ \e^s \ \sum_{i=0}^{m}
\binomial{n}{s+i}\ \mdeg_i(M) \ \sigma^{m-i} \e^i,
\end{equation*}
with equality for $\e$ small enough. 
\end{theorem}

\begin{Proof}
In light of the tube formula, Theorem~\ref{thm:tubeineq}, we aim to bound the integrals of absolute curvature $|K_i|$ of $M$. 
By the Crofton's formula~(\ref{eq:croftonabs}) and the degree bound, Lemma~\ref{le:curvedegree}, we have
{\small
\begin{align*}
  \frac{1}{s+i} |K_i|(M)&\leq \frac{2}{s+i}\binomgr{n}{s+i}\int_{L\in\mathcal{E}_{s+i}^{n}} |K_i|(M\cap L)d\lambda_{s+i}^n\\
&\leq\frac{2 \omega_n}{(s+i)\omega_{s+i}\omega_{m-i}}\binomial{n}{s+i}\Oh_{s+i-1}\int_{L\in\mathcal{E}_{s+i}^{n}}\mdeg(M;L)d\lambda_{s+i}^n.
\end{align*}}

Since $M\subset B^p(p,\sigma)$, we only need to worry about those $L$ that intersect this ball. By our normalization,
\begin{equation*}
\lambda_{s+i}^n(\{L\in \mathcal{E}_{s+i}^n\mid L\cap B^n(p,\sigma)\neq \emptyset\})=\sigma^{m-i}\omega_{m-i},
\end{equation*}
and we can bound the right-most integral above as
\begin{equation*}
\int_{L\in\mathcal{E}_{s+i}^{n}}\mdeg(M;L)d\lambda_{s+i}^n\leq \sigma^{m-i}\omega_{m-i}\cdot \mdeg_i(M).
\end{equation*}
Plugging these bounds into the tube formula~(\ref{eq:weyl}) and simplifying the constants, the claim follows.
\end{Proof}

\subsection{Complete intersections}
Let $f_1,\dots,f_s\in \R[X_1,\dots,X_n]$ be polynomials such that their common zero set $V$ is a complete
intersection, i.e., for every $p\in V$ the gradients $\nabla f_1,\dots,\nabla f_s$ are linearly independent. 
The gradients determine an orientation of $V$.

\begin{lemma}\label{le:compint}
Let $V$ be a complete intersection defined as the zero-set of polynomials $f_1,\dots,f_s$ of degree at most $D$. Then
the degree of the generalized Gauss map $\gamma \colon S(NV)\rightarrow  S^{n-1}$ is bounded by
\begin{equation*}
  \mdeg \gamma \leq (2D)^n.
\end{equation*}
\end{lemma}

\begin{Proof}
We assume $V$ is compact, the general case can be handled with some care.
Let $f=\sum_{i=1}^s f_i^2$, so that in particular, $Z(f)=Z(f_1,\dots,f_s)$. Let $\delta>0$ be such that $\delta$ is a regular value of $f\colon \R^n\rightarrow \R$ and set $f_\delta=f-\delta$.
Then $V_\delta=Z(f_\delta)$ is a hypersurface with associated Gauss map $\gamma_\delta(x)=\nabla f_\delta(x)/\|\nabla f_\delta(x)\|$. By a standard argument using B\'ezout's Theorem (c.f.,~\cite{miln:64}), the degree
of $\gamma_\delta$ is bounded by $(2D)^n$. We next argue that this bound also applies to the cardinality of $\gamma^{-1}(v)$. 

In fact, we can find a regular value of $f$, $\delta>0$, such that for 
all $(p_i,v_i)\in \gamma^{-1}(v)$ and disjoint some neighborhoods $U_i$ of $p_i$ in $\R^n$, there exist $q_i\in U_i$ such that $f(q_i)=\delta$ and $\gamma_\delta(q_i)=v$.
It follows that the number of points in the preimage $\gamma^{-1}(v)$ is also bounded by $(2D)^n$.
\end{Proof}

Now everything is in place for the proof of the main bound.

\Proofof{Theorem~\ref{th:main}}
Set $M':=V\cap B^n(p,\sigma+\e)$. For almost all
$\sigma$, $M'$ will be a smooth compact Riemannian manifold with
smooth $(m-1)$-dimensional boundary $\partial M'$ (Lemma \ref{le:transverse}). Moreover,
\begin{equation*}
T(V,\e)\cap B^n(p,\sigma)\subseteq T(M'\backslash\partial M',\e)\cup T(\partial M',\e). 
\end{equation*}
Note that this inclusion does not hold if we had defined $M'$ by intersecting $V$ with $B^n(p,\sigma)$, as $V$ need not intersect that ball at all.
We can then apply Theorem~\ref{thm:degbound} to $M'\backslash \partial M'$ and to $\partial M'$. 

Since $\mdeg_i(M'\backslash \partial M')\leq \mdeg_i(V)$, it remains to bound $\mdeg_i(V)$. To bound the degree of $V\cap L$, after a change of coordinates we can assume 
that $L$ is given by $x_{s+i+1}=0,\dots,x_{n}=0$. The $f_i$ can therefore be seen as polynomials in $s+i$ variables, denoted by $\overline{x}$. 
The claim now follows from Lemma~\ref{le:compint}.

The boundary $\partial M'$ is defined by the same set of polynomials $f_1,\dots,f_s$ as $V$, with the additional constraint of lying on the sphere $\sum_i x_i^2=(\sigma+\e)^2$.
We can therefore apply the same degree bounds, with the exponents increased by one, to this set.

Note that if $V$ is homogeneous, we can define $M'$ by intersecting with $B^n(p,\sigma)$ rather than $B^n(p,\sigma+\e)$.
We also have $T(V,\e)\cap B^n(p,\sigma)\subseteq T(M'\backslash \partial M',\e)$, which accounts for the factor of $2$ instead of $4$ 
and the simpler form in the second equation in Theorem~\ref{th:main}.

Dividing the resulting expressions by $\vol\ B^n(p,\sigma)=\omega_n\sigma^n$ gives the desired bounds.
\endProofof


\section*{Appendix}
In this appendix we give a proof of the tube formula Theorem~\ref{thm:tubeineq}.

\Proofof{Theorem~\ref{thm:tubeineq}}
We prove the first inequality and point out on the way how the equality for small $\e$ is obtained.
We restrict to compact manifolds without boundary $M$, extending the argument the slightly more general case in the statement of the theorem causes no problem.

Consider the surjective map
\begin{align*}
  f\colon S(\normal M)\times [0,\e] & \rightarrow T(M,\e)\subseteq
  \R^n\\
(p,v,t) &\mapsto p+tv
\end{align*}
of compact manifolds.
For $(p,v)\in S(NM)$ the critical radius is defined as
\begin{equation*}
  \rho_M(p,v)=\sup \{ t \mid \mathrm{dist}(p+tv,M)=t\}, 
\end{equation*}
and set $\rho_M=\mathrm{inf}_{(p,v)\in S(NM)}\rho_M(p,v)$. 
The map $f$ is injective if $\e \leq \rho_M$.

By Sard's Theorem the set of critical values of $f$ has Lebesgue measure zero and
the fibers of $f$ at regular values are finite and locally constant~\cite[\S 1]{miln:97}. 
Given the natural volume form $\omega_{\R^n}$ on $\R^n$ we thus have, by~(\ref{eq:coarea}),
\begin{equation}\label{eq:trans}
  \vol\ T(M,\e)\leq \int_{p\in T(M,\e)}\# f^{-1}(p)
  \ \omega_{\R^n}=\int_{S(\normal M)\times (0,\e)}f^*\omega_{\R^n},
\end{equation}
with equality if $\e\leq \rho_M$. Recall that we are dealing with unsigned forms.

The problem reduces to evaluating the right-hand side. We claim that
\begin{equation}\label{eq:toshow}
  f^*\omega_{\R^n}=t^{s-1}|\det(\mathrm{Id}-tS(v))\ \omega_{S(\normal M)}\wedge dt|.
\end{equation}
Assuming this to hold for the moment, the claimed inequality for the
volume of tubes follows by integrating
\begin{align*}
  \int_{S(\normal M)\times (0,\e)}f^*\omega_{\R^n} &=\int_{S(\normal M)} \left(\int_{0}^{\e}
  t^{s-1}|\det(\mathrm{Id}-tS(v))| \ dt \right) \omega_{S(\normal M)}\\
&\leq \int_{S(\normal M)}
  \left(\int_{0}^{\e}t^{s-1}\sum_{i=0}^{m}t^i|\psi_i(v)|
  \ dt\right)\omega_{S(\normal M)}\\
&=\sum_{i=0}^m \left(\int_{0}^{\e}t^{s-1+i}\ dt\right)
  \left(\int_{S(\normal M)}|\psi_i(v)|\ \omega_{S(\normal M)}\right) \\
  &= \sum_{i=0}^{m}
  \frac{1}{s+i} \ \e^{s+i} \ |K_i|(M).\\ 
\end{align*}
It therefore remains to prove (\ref{eq:toshow}). Note that is $\e<\rho_M$, then the map $f$ is injective and, with the right choice of orientation,
the determinant $\det(\mathrm{Id}-tS(v))$ is always positive. We can therefore omit the absolute value and obtain an equality with the integrals of curvature.

Let $(x^1,\dots,x^m)\colon U\rightarrow \R^m$ be orthonormal
coordinates on $U=U'\cap M$. Let $(E_1,\dots,E_n)$ be an
orthonormal frame field on $U'$ such that
$E_i:=~\frac{\partial}{\partial x^i}$ on $U\subseteq M$ for $1\leq i\leq m$.
Set $\omega_M:=E_1^*\wedge\cdots \wedge E_m^*$ and 
$\omega_N:=E_{m+1}^*\wedge \cdots \wedge E_n^*$ ($E_i^*$ denoting the
dual of $E_i$). We then have $\omega_{\R^n}=\omega_M\wedge \omega_N$, and
for the restriction to $M$, $\omega_M|_{TM}=dx^1\wedge \cdots \wedge dx^m$.

The frame field also gives a local trivialization of the sphere bundle
\begin{align*}
  U\times S^{s-1}&\rightarrow S(\normal M)\\
  (p,u)&\mapsto \left(p,\sum_{i=1}^{s}u^i E_{m+i}(p)\right).
\end{align*}
An orthonormal coordinate system $y^1,\dots,y^{s-1}$ for $S^{s-1}$
then gives rise to orthonormal coordinates
$(x^1,\dots,x^m,y^1,\dots,y^{s-1},t)$ on $S(\normal M)\times
(0,\e)$. Setting $dx=dx^1\wedge \cdots \wedge dx^{m}$ and $dy=dy^1\wedge
\cdots \wedge dy^{s-1}$ we have
\begin{equation}\label{eq:wedge1}
  \omega_{S(\normal M)}\wedge dt=dx\wedge dy\wedge dt.
\end{equation}
Let $\phi(p,v,t)$ be such that
$f^*\omega_{\R^n}=\phi(p,v,t) \ \omega_{S(\normal M)}\wedge dt$ as differential form.
By Equation (\ref{eq:wedge1}) we obtain
\begin{align*}
  \phi(p,v,t)&=f^*\omega_{\R^n}\left(
\frac{\partial}{\partial x^1},\dots,\frac{\partial}{\partial
  x^m},\frac{\partial}{\partial y^1},\dots,\frac{\partial}{\partial
  y^{s-1}},\frac{\partial}{\partial t}\right)\\
&=\omega_{\R^n}\left(f_*\frac{\partial}{\partial x^1},\dots,f_*\frac{\partial}{\partial
    x^m},f_*\frac{\partial}{\partial
    y^1},\dots,f_*\frac{\partial}{\partial
    y^{s-1}},f_*\frac{\partial}{\partial t}\right).
\end{align*}

We next observe that, using the definition of $f$,
\begin{align*}
f_*\frac{\partial}{\partial x^i}&=\frac{\partial}{\partial x^i}p+t\sum_{\ell=1}^s u^{\ell}\frac{\partial}{\partial x^i}E_{m+\ell}(p),\\
f_*\frac{\partial}{\partial y^j}&=t\sum_{\ell=1}^s \frac{\partial}{\partial y^j}u^{\ell} \ E_{m+\ell}(p),\\
f_*\frac{\partial}{\partial t}&=v.
\end{align*}
In particular,
$f_*(T_vS^{s-1}\times T_t\R)\subseteq N_pM=(T_pM)^{\perp}$, so that
\begin{equation*}
\phi(p,v,t)=\omega_{M}\left(\frac{\partial}{\partial x^1} f,\dots,\frac{\partial}{\partial x^m} f\right) \cdot
  \omega_{N}\left(\frac{\partial}{\partial y^1} f,\dots,\frac{\partial}{\partial y^{s-1}}
  f,\frac{\partial}{\partial t} f\right).
\end{equation*}

A straight-forward calculation shows that
\begin{equation*}
  \left\langle\frac{\partial}{\partial x^i}f,E_j\right\rangle=\left\langle E_i+t\frac{\partial}{\partial x^i}Z,E_j\right\rangle=\delta_{ij}-tS_{ij}(v),
\end{equation*}
where $Z$ is a normal vector field with $Z(p)=v$.
From this it follows that 
\begin{equation*}
\omega_{M}\left(\frac{\partial}{\partial x^1}
f,\dots,\frac{\partial}{\partial x^m} f\right)=\det
(\mathrm{Id}-tS(v)).
\end{equation*}
Similarly one obtains 
\begin{equation*}
\omega_{N}\left(\frac{\partial}{\partial
    y^1} f,\dots,\frac{\partial}{\partial y^{s-1}}
  f,\frac{\partial}{\partial t} f\right)=t^{s-1}.
\end{equation*}
This completes the proof of the claimed inequality.
%
\endProofof

\bibliographystyle{plain}
\bibliography{tubes-v2}

\end{document}